# Elliptic curves of high rank over function fields

Jasper Scholten

### Abstract

By the Mordell-Weil theorem the group of $\mathbf{Q}(z)$-rational points of an elliptic curve is finitely generated. It is not known whether the rank of this group can get arbitrary large as the curve varies. Mestre and Nagao have constructed examples of elliptic curves $E$ with rank at least 13.

In this paper a method is explained for finding a 14th independent point on $E$, which is defined over $k(z)$, with $[k : \mathbf{Q}] = 2$. The method is applied to Nagao's curve. For this curve one has $k = \mathbf{Q}(\sqrt{-3})$.

The curves $E$ and 13 of the 14 independent points are already defined over a smaller field $k(t)$, with $[k(z) : k(t)] = 2$. Again for Nagao's curve it is proved that the rank of $E(\bar{\mathbf{Q}}(t))$ is exactly 13, and that rank $E(\mathbf{Q}(t))$ is exactly 12.

## 1 Mestre's construction

First a method due to Mestre [2] for constructing elliptic curves with high rank is described.

Let $k$ be any field with char $k \neq 2$. Choose $2n$ elements $a_1, \ldots, a_{2n} \in k$. We are going to construct a plane curve $C$ such that the points with $x$-coordinate $a_i$ are $k$-rational. To do so, set $p(x) = \prod_{i=1}^{2n}(x - a_i)$. It is easily shown that there exist polynomials $q$ and $r$ in $k[x]$ with $\deg r \leq n - 1$ such that $p = q^2 - r$. Define $C$ by the equation $y^2 = r(x)$. Clearly $C$ contains the points $(a_i, \pm q(a_i))$.

For $n = 5$ almost all choices for the $a_i$ give that $\deg r = 4$ and that $C$ is a curve of genus 1 with 10 points of the form $(a_i, \pm q(a_i))$. If $C$ is made into an elliptic curve by choosing one of these points as the zero point then the other points generate a group of rank 9 (generically).

Mestre constructs an elliptic curve of rank 11 over $\mathbf{Q}(t)$ by taking $n = 6$ and $a_i = b_i + t$ for $i = 1, \ldots, 6$, and $a_i = b_{i-6} - t$ for $i = 7, \ldots, 12$. Now the $x^5$-coefficient of $r$ is of the form $s \cdot t^2$ with $s \in \mathbf{Q}[b_1, \ldots, b_6]$. It is not very difficult to find $b_i \in \mathbf{Q}$ such that $s(b_1, \ldots, b_6) = 0$. (One can for example choose $b_1$ upto $b_5$ at random and hope that there exists a $b_6 \in \mathbf{Q}$ with $s(b_1, \ldots, b_6) = 0$. Trying this often enough will almost surely give the desired $b_i$'s.) In Mestre's example we have $b_1 = -17$, $b_2 = -16$, $b_3 = 10$, $b_4 = 11$, $b_5 = 14$, $b_6 = 17$. If one chooses one of the points $(a_i, \pm q(a_i))$ as zero point then the other points generate a group of rank 11. In fact, it is easy to see that they cannot generate a group of larger rank: the divisors of the functions $(x - a_i)/(x - a_j)$ give relations $(a_i, q(a_i)) + (a_i, -q(a_i)) = (a_j, q(a_j)) + (a_j, -q(a_j))$, and the divisor of the function $(q(x) - y)/(q(x) + y)$ gives the relation

$$\sum_{i=1}^{12}(a_i, q(a_i)) = \sum_{i=1}^{12}(a_i, -q(a_i)).$$

The curve Mestre has constructed this way is

$$\begin{aligned}
y^2 = {}& \left(213040 + 429\,t^2\right) x^4 + \left(-4956000 - 5434\,t^2\right) x^3 + \\
& \left(4888624 - 2451\,t^2 - 858\,t^4\right) x^2 + \left(539121408 - 3637984\,t^2 + 5434\,t^4\right) x \\
& - 3035397056 + 53200096\,t^2 - 268748\,t^4 + 429\,t^6.
\end{aligned}$$

(Here the polynomial $r$ is replaced by $4r/81t^2$ in order to obtain a slightly simpler equation.)

From this curve one can construct a new curve with 12 independent points (see [3]). To do so, parametrise the conic $u^2 = 213040 + 429\,t^2$. This can be done as follows:

$$t = \frac{6\,z^2 - 956\,z + 2574}{z^2 - 429}, \quad u = \frac{205062 + 478\,z^2 - 5148\,z}{z^2 - 429}.$$

Now the leading coefficient of $r$ is a square in $\mathbf{Q}(z)$, so the points at infinity have become $\mathbf{Q}(z)$-rational. These points together with the points we already had generate a group of rank 12.

## 2 Nagao's 13th point

In [5] Nagao constructs an elliptic curve of rank at least 13 over $\mathbf{Q}(t)$. He uses Mestre's method to construct a curve of rank at least 12, but he takes a different choice for the $b_i$. He takes $b_1 = 148$, $b_2 = 116$, $b_3 = 104$, $b_4 = 57$, $b_5 = 25$, $b_6 = 0$, yielding the curve

$$\begin{aligned}
y^2 = {}& \left(330112972800 + 14017536\,t^2\right) x^4 + \\
& \left(-4205260800\,t^2 - 99527168931840\right) x^3 + \\
& \left(10445363957523456 + 617005209600\,t^2 - 28035072\,t^4\right) x^2 + \\
& \left(-443009190070886400 - 46063725926400\,t^2 + 4205260800\,t^4\right) x + \\
& 6473450277365760000 + 1388825681338368\,t^2 - 247954636800\,t^4 + \\
& 14017536\,t^6.
\end{aligned} \qquad (1)$$

With this choice he finds the extra, independent point $((t+703)/15, (-224\,t^3 - 844\,t^2 + 900484\,t + 2161725)/75)$. Here the point at infinity is rational over $\mathbf{Q}(z)$, with $t = \frac{23550 - z^2}{2\,z}$.

Mestre has found a 2-parameter family of curves over $\mathbf{Q}(t)$ with an independent 13th point. Nagao's curve is in this family. For the exact equations of this family we refer to [4].

## 3  Elliptic surfaces

Let us give a brief review on the theory of elliptic surfaces. Everything we need is in [7]. For us an elliptic surface will mean the following: a projective regular surface $S$ defined over an algebraically closed field $k$ together with a map $f : S \longrightarrow C$ to a regular projective curve $C$ such that the generic fibre of $f$ is a curve of genus 1 over the function field $k(C)$ of $C$, and such that

1. $S$ is relatively minimal, i.e. the fibres of $f$ contain no exceptional curves of the first kind,

2. $S$ does not split, i.e. $S$ is not $C$-birationally equivalent to a constant family $E_0 \times C$, and

3. $\phi$ has a section $O$, i.e. a map $O : C \longrightarrow S$ with $f \circ O = \operatorname{Id}_C$.

If $E$ is an elliptic curve over $k(C)$ then there exists an elliptic surface $f : S \longrightarrow C$ with generic fibre $E$, and it is unique upto $C$-isomorphism. Moreover, the points in $E(k(C))$ correspond bijectively to the sections of $f$. We take $O$ to be the section corresponding to the zero point of $E$. For $P \in E(k(C))$ we define $(P)$ to be the image of the corresponding section. So $(P)$ is a divisor on $S$.

The Néron-Severi group $\mathrm{NS}(S)$ of $S$ is the group of divisors of $S$ modulo algebraic equivalence. On $\mathrm{NS}(S)$ we have a non-degenerate integral bilineair pairing, the intersection pairing. Let us denote it by $(\,.\,,\,.\,)$.

If we let $m_v$ denote the number of components of $f^{-1}(v)$, for $v \in C$, then the following relation holds:

$$\operatorname{rank} E(k(C)) = \operatorname{rank} \mathrm{NS}(S) - 2 - \sum_{v \in C}(m_v - 1). \tag{2}$$

The Néron-Tate height $h$ on the Mordell-Weil group $E(k(C))$ can be expressed in terms of certain intersection numbers. We have

$$h(P) = 2((P), (O)) - 2((O), (O)) - \sum_{v \in C} \lambda_v(P),$$

where $\lambda_v(P)$ is a rational number depending on which component of $f^{-1}(v)$ is being hit by $(P)$. It is equal to 0 if $(P)$ hits the same component as $(O)$, so that is almost everywhere since almost all fibres are irreducible. The only other case where we will need to know the value of $\lambda_v$ is when $f^{-1}(v)$ is of Kodaira type $I_2$. In this case $\lambda_v(P) = \frac{1}{2}$ if $(P)$ and $(O)$ hit different components.

In general the rank of the Néron-Severi group is not easy to compute. However, if the surface $S$ is a rational surface then we have $\text{rank}\,\text{NS}(S) = 10$. So in this case computing the rank of the Mordell-Weil group is equivalent to computing the number of components in reducible fibres. The latter can be done very easily with Tate's algorithm.

It turns out that if $S$ is rational then not only the rank of $E(k(C))$ is determined by the reducible fiber types. The lattice structure on $E(k(C))$ induced by the Néron-Tate height is also determined by the reducible fibre configuration, except for a few cases where there are 2 possible lattice structures. See [6]. A particular example which will be used in this sequel is: if $S$ is rational, and $f : S \to C$ has exactly one reducible fiber, of type $I_2$, then the lattice on $E(k(C))$ induced by the Néron-Tate height is the dual lattice $E_7^*$ of the root lattice $E_7$.

Checking whether an elliptic surface $S$ with generic fibre $E/k(C)$ is rational is very easy. A necessary condition for $S \to C$ to be rational is that the curve $C$ is rational, so $k(C) = k(t)$. Now $S$ is rational if and only if there is a Weierstrass equation for $E/k(t)$ such that the Weierstrass coefficients $a_i$ of this equation are polynomials of degree less than or equal to $i$.

## 4 A 14th point

In this section a method for finding an extra independent point on Mestre's rank 12 curves is discussed. This point will be defined over $k(t)$, with $k$ a quadratic extension of $\mathbf{Q}$. If one starts with a curve of rank at least 13, as descibed in section 2, then a curve of rank at least 14 is obtained.

I will explain the method by looking at Nagao's curve (1). Note that all coefficients of the equation defining the curve are even polynomials in $t$, so this curve is defined over $\mathbf{Q}(u)$, with $u = t^2$. Then a minimal Weierstrass equation for this curve will also be defined over $\mathbf{Q}(u)$ (See e.g. [1], chapter 20). Using $(t, 2544297600 - 87059232\,t + 836160\,t^2)$ as zero point one finds this equation:

$$
\begin{aligned}
y^2 &= x^3 + (-432\,u^4 - 4435200\,u^3 + 38353513056\,u^2 - \\
&\quad 18899197014000\,u - 340079781902569707)x + \\
&\quad 3456\,u^6 + 53222400\,u^5 - 870054636672\,u^4 + \\
&\quad 5893342291009600\,u^3 - 18532375351306853196\,u^2 + \\
&\quad 7556017995191414902800\,u + 7630232732696966177153 1494.
\end{aligned} \tag{3}
$$

The Weierstrass coefficients $a_i$ have degree less than or equal to $i$, so the corresponding elliptic surface $S$ is rational. It follows from Tate's algorithm that $S \to \mathbf{P}^1$ has exactly one reducible fibre. This reducible fibre is of Kodaira type $I_2$ and lies above the point $u = \infty \in \mathbf{P}^1$. Hence the Mordell-Weil lattice $E(\bar{\mathbf{Q}}(u))$ is equal to $E_7^*$. Here $\bar{\mathbf{Q}}$ denotes the algebraic closure of $\mathbf{Q}$.

Consider the norm map

$$\begin{aligned} N : E(\mathbf{Q}(t)) &\longrightarrow E(\mathbf{Q}(u)) \\ P &\longmapsto P + \sigma(P), \end{aligned}$$

where $\sigma$ is the non-trivial element in $\mathrm{Gal}(\mathbf{Q}(t)/\mathbf{Q}(u))$. Having explicit generators of a subgroup $W \subset E(\mathbf{Q}(t))$ of rank 12, one can compute generators of the image $N(W) \subset E(\mathbf{Q}(u))$. Using the height pairing on $E(\mathbf{Q}(u))$ one can also compute the rank of $N(W)$. This rank is equal to 6. So there must be an extra point in $E(\bar{\mathbf{Q}}(u))$, independent of $N(W)$. Let $Q \in E(\bar{\mathbf{Q}}(u))$ be such an extra point.

**Lemma 1** *The point $Q$ is in $E(k(u))$ for some quadratic extension $k$ of $\mathbf{Q}$.*

*Proof.* Let $T$ denote the primitive closure of $N(W)$ in $E(\bar{\mathbf{Q}}(u))$, i.e.

$$T = \{P \in E(\bar{\mathbf{Q}}(u)) \mid nP \in N(W) \text{ for some } n \geq 1\}.$$

Let $P \in T$ such that $nP \in N(W)$, and let $\sigma \in \mathrm{Gal}(\bar{\mathbf{Q}}/\mathbf{Q})$. Then

$$0 = nP - \sigma(nP) = n(P - \sigma(P)) \in E(\bar{\mathbf{Q}}(u)) \cong \mathbf{Z}^7.$$

So $P$ is $\sigma$-invariant, and $\mathrm{Gal}(\bar{\mathbf{Q}}/\mathbf{Q})$ acts trivial on $T$. Now it suffices to show that the representation

$$\mathrm{Gal}(\bar{\mathbf{Q}}/\mathbf{Q}) \longrightarrow \mathrm{Aut}(E(\bar{\mathbf{Q}}(u))/T)$$

has a kernel of index 1 or 2. But that follows from $E(\bar{\mathbf{Q}}(u))/T \cong \mathbf{Z}$ so $\mathrm{Aut}(E(\bar{\mathbf{Q}}(u))/T) = \{\pm \mathrm{Id}\}$. □

The lattice $E_7^*$ is generated by vectors of length $3/2$, so we may assume that $Q$ has length $3/2$. From

$$h(Q) = 2((Q),(O)) - 2((O),(O)) - \sum_{v \in C} \lambda_v(Q)$$

and from the fact that for a rational elliptic surface we have $((O),(O)) = -1$ (see [7]) we see that $(Q)$ must be disjoint from the zero section, and it must hit the $I_2$-fibre in the non-zero component. This information can be used to find $Q$. From $((Q),(O)) = 0$ it follows that $Q = (X(u), Y(u))$ with $X(u)$ and $Y(u)$ in $k[u]$. In order to examine the fiber at infinity we make a coordinate change

$$u' = \frac{1}{u}, \quad x' = x\,u'^2, \quad y' = y\,u'^3.$$

This yields a new Weierstrass equation which at $u' = 0$ is

$$y'^2 = (x' - 12)^2 (x' + 24).$$

The point $(x', y') = (12, 0)$ is a node. Since $Q$ hits the non-zero component of the fibre at $u = \infty$ we have that $(u'^2 X(\frac{1}{u'}),\ u'^3 Y(\frac{1}{u'}))$ passes through this node. So

$$\left.(u'^2 X(\frac{1}{u'}),\ u'^3 Y(\frac{1}{u'}))\right|_{u'=0} = (12, 0),$$

and thus $X(u)$ has degree 2 and leading term $12\,u^2$, and $Y(u)$ has degree at most 2. Write $X(u) = 12\,u^2 + a\,u + b$ and $Y(u) = c\,u^2 + d\,u + e$. Substituting this in the equation (3) yields a system of equations with unknowns $a, b, c, d$ and $e$. The computer algebra system Maple was able to solve this system. One of the solutions is

$$\begin{aligned}
X(u) &= 12\,u^2 + \left(129032 + 27432\,\sqrt{-3}\right)u - 168316272\,\sqrt{-3} - 757109813, \\
Y(u) &= \left(164592 + 404592\,\sqrt{-3}\right)u^2 + \left(1518893856 - 3588498678\,\sqrt{-3}\right)u \\
&\quad - 15516067218048 + 7021689895536\,\sqrt{-3}.
\end{aligned}$$

This point can be taken as $Q$.

**Theorem 1** *The point $Q$ together with Nagao's 13 points generate a group of rank 14 over $\mathbf{Q}(\sqrt{-3}, z)$ (with $z$ defined as in section 2).*

*Proof.* Let $\sigma$ denote the non-trivial element in $\mathrm{Gal}(\mathbf{Q}(\sqrt{-3}, z)/\mathbf{Q}(z))$. We have

$$0 \neq Q - \sigma(Q) \in E(\bar{\mathbf{Q}}(u)) \cong \mathbf{Z}^7,$$

so $Q - \sigma(Q)$ has infinite order. It follows that $nQ \notin E(\mathbf{Q}(\sqrt{-3}, z))^{\sigma = \mathrm{Id}} = E(\mathbf{Q}(z))$ for all $n \geq 1$. So $Q$ is independent of Nagao's 13 points. $\square$

## 5 A curve of rank exactly 13

In the previous section a 13th independent point was constructed on Nagao's curve (1). In this section it is shown that there is no 14th independent point on this curve.

**Theorem 2** *Let $E$ be Nagao's curve* (1). *Then $E(\bar{\mathbf{Q}}(t))$ has rank exactly 13.*

*Proof.* Let $S \longrightarrow \mathbf{P}^1$ be the corresponding elliptic surface. From Tate's algorithm it follows that $S$ has only 1 reducible fibre, at $t = \infty$. This fibre is of Kodaira type $I_4$. So in view of (2) it suffices to show that rank $\mathrm{NS}(S) = 18$.

If the surface $S$ has good reduction modulo a prime $p$ then for the reduced surface $\bar{S}$ one has rank $\mathrm{NS}(\bar{S}) \geq$ rank $\mathrm{NS}(S)$. So the theorem is proved if there is a prime $p$ of good reduction such that rank $\mathrm{NS}(\bar{S}) = 18$.

Let $H^i(\bar{S}_{\bar{\mathbf{F}}_p}, \mathbf{Q}_l)$ be the $i$th $l$-adic cohomology group, for some prime $l \neq p$. Denote the Frobenius endomorphism on $H^i(S_{\bar{\mathbf{F}}_p}, \mathbf{Q}_l)$ by $F_p$. From [8] it follows that the number of eigenvalues of $F_p : H^2(S_{\bar{\mathbf{F}}_p}, \mathbf{Q}_l) \to H^2(S_{\bar{\mathbf{F}}_p}, \mathbf{Q}_l)$ of the form $p\zeta$, with $\zeta$ a root of unity, is an upper bound for rank $\mathrm{NS}(\bar{S})$.

These eigenvalues can be computed using the Lefschetz trace formula

$$\#\bar{S}(\mathbf{F}_{p^n}) = \sum_{i=0}^{4}(-1)^i \mathrm{trace}(F_p^n | H^i(\bar{S}_{\bar{\mathbf{F}}_p}, \mathbf{Q}_l)). \tag{4}$$

For a surface we have $\mathrm{trace}(F_q^n|H^0) = 1$ and $\mathrm{trace}(F_p^n|H^4) = p^{2n}$. For a fibred non-isotrivial surface with base $\mathbf{P}^1$ and a section the $H^1$ and $H^3$ vanish. A quadratic base change of a rational semi-stable elliptic surface yields a $K3$ surface, and hence we have dim $H^2 = 22$. The determinant $\det(F_p|H^2)$ equals $\pm p^{22}$.

We already know 18 eigenvalues of $F_p$. There are 4 eigenvalues $p$ coming from the fibre components (this follows from the fact that $E/\mathbf{Q}(t)$ has split multiplicative reduction at $t = \infty$, so the components of the $I_4$ fibre are $\mathbf{Q}$-rational), 13 eigenvalues $p$ coming from the sections defined over $\mathbf{F}_p$ (including the zero section), and 1 eigenvalue $\left(\frac{-3}{p}\right) p$ coming from the section defined over $\mathbf{F}_p(\sqrt{-3})$. All this information enables us to compute the remaining 4 eigenvalues if we know the sign of $\det(F_p|H^2)$, and $\#\bar{S}(\mathbf{F}_{p^n})$ for $n = 1$, 2 and 3.

In order to prove the theorem it is not necessary to compute the 4 remaining eigenvalues. It suffices to know that they are not of the form $p\zeta$. With some luck the latter can be shown without having to perform the time consuming computation of $\#\bar{S}(\mathbf{F}_{p^3})$. For this one reasons as follows:

The 4 eigenvalues are roots of a quartic polynomial in $\mathbf{Q}[X]$. So if one of them is of the form $p\zeta$, then $\zeta$ is a primitive 1st, 2nd, 3rd, 4th, 5th, 6th, 8th, 10th or 12th root of unity. Compute $\#\bar{S}(\mathbf{F}_p)$ and $\#\bar{S}(\mathbf{F}_{p^2})$, and assume that one of the 4 unknown eigenvalues is of the form $p\zeta$. For each choice of $\zeta$, and each choice of $\mathrm{sgn}(\det(F_p|H^2))$ the remaining 3 eigenvalues can be computed using (4). These 3 eigenvalues should contain the conjugates of $p\zeta$, and they should have absolute values $p$ according to Deligne's theorem. If, for all 18 choices of $\zeta$ and $\mathrm{sgn}(\det(F_p|H^2))$, this is not the case then the assumption that one eigenvalue is of the form $p\zeta$ was false.

For $p$ we can take primes with the property that $S$ and $\bar{S}$ have the same singular fibre configuration. The smallest prime for which this holds is $p = 53$. Calculations at this prime show that $\#\bar{S}(\mathbf{F}_{53}) = 3593$ and $\#\bar{S}(\mathbf{F}_{53^2}) = 7945269$. The assumption that 53 is one of the eigenvalues, and that $\det(F_{53}|H^2) = 53^{22}$ implies that the other 3 eigenvalues are roots of

$$X^3 + 118X^2 + 6254X + 148877 = (X + 53)(X^2 + 65X + 2809). \qquad (5)$$

This is consistent with Deligne's theorem. So we don't arrive at the contradiction hoped for. In fact, these eigenvalues are correct, as can be shown by computing $\#\bar{S}(\mathbf{F}_{53^3})$. It is also interesting to remark that from the Tate conjectures (which are proved for $K3$ surfaces) it now follows that $\operatorname{rank}\operatorname{NS}(\bar{S}) = 20$.

The next prime where the singular fibre configuration does not change is $p = 71$. One has $\#\bar{S}(\mathbf{F}_{71}) = 6096$ and $\#\bar{S}(\mathbf{F}_{71^2}) = 25498920$. It turns out that the assumption that one of the 4 remaining eigenvalues is of the form $p\zeta$ leads to a contradiction as mentioned above, for every possible choice of $\zeta$. So all 4 eigenvalues are not of that form, and $18 \leq \operatorname{rank}\operatorname{NS}(S) \leq \operatorname{rank}\operatorname{NS}(\bar{S}) \leq 18$. $\square$

**Corollary 1** *The group $E(\mathbf{Q}(t))$ has rank exactly 12.*

*Proof.* In the proof of theorem 1 it is shown that $nQ \notin E(\mathbf{Q}(t))$ for all $n \geq 1$. This shows that $\operatorname{rank} E(\mathbf{Q}(t)) \neq 13$, and the corollary follows. $\square$

One might try to find an upper bound for $\operatorname{rank} E(\bar{\mathbf{Q}}(z))$ in the same way as it is done here for $\operatorname{rank} E(\bar{\mathbf{Q}}(t))$. But the computations become too big. The elliptic surface with generic fibre $E/\mathbf{Q}(z)$ has an $H^2$ of dimension 46, and only 26 eigenvalues of Frobenius are known in advance. So one would have to count points over finite fields with cardinality upto $p^{18}$.

Jasper Scholten
Vakgroep Wiskunde
Rijksuniversiteit Groningen
postbus 800
9700 AV Groningen
the Netherlands
e-mail: jasper@math.rug.nl